\documentclass[a4paper,leqno,12pt]{article}
\usepackage{amsmath,amsthm,amsfonts,amssymb}

\makeatletter
\def\@seccntformat#1{\csname the#1\endcsname.\quad}
\renewcommand\section{\@startsection {section}{1}{\z@}%
                                   {-3.5ex \@plus -1ex \@minus -.2ex}%
                                   {2.3ex \@plus.2ex}%
                                   {\normalfont\bf\center }}
\renewcommand\subsection{\@startsection {subsection}{1}{\z@}%
                                   {-3.5ex \@plus -1ex \@minus -.2ex}%
                                   {2.3ex \@plus.2ex}%
                                   {\normalfont\bf}}
\makeatother

\title{\large\bf Brownian sheet and reflectionless potentials}
\author{\normalsize Setsuo Taniguchi
\thanks{Research supported partially by Grant-in-Aid for
     Scientific Research (A)(1) 14204010}
\\
{\normalsize\it  Faculty of Mathematics, Kyushu University,
       Fukuoka 812-8581, Japan}}
\date{\normalsize July 12, 2005}

\newtheorem{thm}{Theorem}
\newtheorem{cor}{Corollary}
\newtheorem{prop}{Proposition}
\newtheorem{lem}{Lemma}
\theoremstyle{definition}
\newtheorem{rem}{Remark}

\begin{document}
\maketitle

\begin{abstract}
In this paper, the investigation into stochastic calculus
related with the KdV equation, which was initiated by
S.~Kotani \cite{kotani} and made in succession by N.~Ikeda
and the author \cite{ikeda-taniguchi,taniguchi}, is
continued.
Reflectionless potentials give important examples in the
scattering theory and the study of the KdV equation; they
are expressed concretely by their corresponding scattering
data, and give a rise of solitons of the KdV equation.
N.~Ikeda and the author \cite{ikeda-taniguchi} established a
mapping $\psi$ of a family $\mathcal{G}_0$ of probability
measures on the $1$-dimensional Wiener space to the space
$\Xi_0$ of reflectionless potentials.
The mapping gives a probabilistic expression of
reflectionless potential.
In this paper, it will be shown that $\psi$ is bijective,
and hence $\mathcal{G}_0$ and $\Xi_0$ can be identified.
The space $\Xi_0$ was extended to the one $\Xi$ of
generalized reflectionless potentials, and was used by
V.~Marchenko to investigate the Cauchy problem for the KdV
equation and by S.~Kotani to construct KdV-flows. 
As an application of the identification of $\mathcal{G}_0$
and $\Xi_0$ via $\psi$, taking advantage of the Brownian
sheet, it will be seen that convergences of elements in 
$\mathcal{G}_0$ realizes the extension of $\Xi_0$ to $\Xi$.

\end{abstract}


\section{Introduction}
\label{intro}
Let $\mathcal{W}$ be the space of all $\mathbf{R}$-valued 
continuous functions $w$ on $[0,\infty)$ with $w(0)=0$, and
$\mathcal{B}$ be its Borel $\sigma$-field,
$\mathcal{W}$ being equipped with the topology of uniform 
convergence on compacts.
The coordinate mapping on $\mathcal{W}$ is denoted by 
$X(x)$; $X(x,w)=w(x)$, $w\in \mathcal{W}$, $x\in[0,\infty)$.
Let $\Sigma_0$ be the set of measures on $\mathbf{R}$ of
the form $\sum_{j=1}^n c_j^2 \delta_{p_j}$ for some
$n\in \mathbf{N}$ and $p_j\in\mathbf{R}$, $c_j>0$, 
$1\le j\le n$ with $p_i\ne p_j$ if $i\ne j$,
where $\delta_p$ is the Dirac measure concentrated at $p$.
For $\sigma\in\Sigma_0$, set 
$$
 R_\sigma(x,y)=\int_{\mathbf{R}} 
  \frac{e^{\zeta(x+y)}-e^{\zeta|x-y|}}{2\zeta} \sigma(d\zeta),
$$
and let $P^\sigma$ be the probability measure on
$(\mathcal{W},\mathcal{B})$ so that
$\{X(x)\}_{x\ge0}$ is a centered Gaussian process with
covariance function $R_\sigma$ 
(a construction of $P^\sigma$ will be given in 
Sect.~\ref{sec:2}).
Put
$$
 \mathcal{G}_0=\{P^\sigma|\sigma\in\Sigma_0\}.
$$
N. Ikeda and the author (\cite{ikeda-taniguchi}) showed
that, for each $P^\sigma$, the function 
\begin{equation}\label{eq.i&t}
  \psi(P^\sigma)(x)= 4\Bigl(\frac{d}{dx}\Bigr)^2 \,\log \biggl(
  \int_{\mathcal{W}} \exp\biggl(-\frac12 \int_0^x X(y)^2 dy
  \biggr)dP^\sigma\biggr),\quad x\ge0,
\end{equation}
is well defined and coincides with the restriction of a
reflectionless potential to $[0,\infty)$, and the associated
scattering data was specified in terms of $\sigma$. 
Since reflectionless potentials are real analytic, we may
and will think of $\psi$ as a mapping of $\mathcal{G}_0$ to
the space $\Xi_0$ of reflectionless potentials.
It should be recalled that reflectionless
potentials give a rise of solitons of the KdV equation
(\cite{mar-book,miwa-jimbo-date}). 
A review on these results and the definition of
reflectionless potential will be given in
Sect.~\ref{sec:2}. 
The first aim of this paper is to show that the mapping
$\psi$ of $\mathcal{G}_0$ to $\Xi_0$ is bijective. 
See Theorem~\ref{thm.bij}.
In this sense, the set $\Xi_0$ of analytic future and the
set $\mathcal{G}_0$ of probabilistic future are identified.
Moreover, we shall establish a probabilistic
expression of $u\in\Xi_0$ through $\psi$.
See Corollaries~\ref{cor.neg} and \ref{cor.brown}.

A generalized reflectionless potential $u$ is a limit of a
sequence $\{u_n\}$ of reflectionless potentials $u_n$ 
such that 
$\mbox{\rm Spec}(-(d/dx)^2+u_n)\subset[-\lambda_0,\infty)$, 
$n=1,2,\dots$, for some $\lambda_0>0$ in the topology of
uniform convergence on compacts, where
$\mbox{\rm Spec}(-(d/dx)^2+u_n)$ denotes the spectrum of 
$-(d/dx)^2+u_n$.
The space $\Xi$ of generalized potentials was used by
V.~Marchenko (\cite{marchenko}) to study the Cauchy problem
for the KdV equation, and by S. Kotani (\cite{kotani}) to
construct KdV-flows. 
Let $\Sigma$ be the space of all finite measures on
$\mathbf{R}$ with compact support, and
$$
 \mathcal{G}=\{P^\sigma|\sigma\in\Sigma\}, 
$$
where we have naturally extended the notation $P^\sigma$ to
$\Sigma$.
On account of the identification of $\Xi_0$ and
$\mathcal{G}_0$ stated in the above paragraph, arises 
a natural question if one can describe the relation between
convergences of reflectionless potentials to generalized
ones and convergences of probability measures in
$\mathcal{G}_0$ to those in $\mathcal{G}$.
The second aim of this paper is to answer affirmatively to
this question. 
Namely, we shall study the convergence of $\psi(P^\sigma)$'s
with $P^\sigma$ not only in $\mathcal{G}_0$ but also in
$\mathcal{G}$. 
In particular, the convergence of elements in $\Xi_0$
defining those in $\Xi$ will be realized through the
convergence of elements in $\mathcal{G}_0$ to those in
$\mathcal{G}$.
Moreover, we shall show that the surjectivity of $\psi$ on
$\mathcal{G}_0$ extends to $\mathcal{G}$; every $u\in\Xi$
admits $P^\sigma\in \mathcal{G}$ so that $\psi(P^\sigma)=u$
on $[0,\infty)$.
The expression of such $u$ on $(-\infty,0]$ by $\psi$ and
$\sigma$ will be also given. 
For these, see Theorem~\ref{thm.conv.gen} and
Remark~\ref{rem.gen.rl}.
A key ingredient for the investigation is to realize the
above $P^\sigma$ by using the Brownian sheet and reduce
every estimations to the ones for Wiener integrals
associated with the Brownian sheet.

The organization of the paper is as follows.
In Sect.~\ref{sec:2}, we shall show the bijectivity of
$\psi$ in (\ref{eq.i&t}) after reviewing the result in
\cite{ikeda-taniguchi}.
In the section, a construction of $P^\sigma$ for
$\sigma\in\Sigma_0$ is given.
Sect.~\ref{sec:3} is devoted to introducing compound
Ornstein-Uhlenbeck processes which are indispensable to 
discuss the convergence of $P^\sigma$'s.
The Brownian sheet plays a key role to construct such
processes.
Another realization of $P^\sigma$ with the Brownian sheet
will be also given there.
In the last section, we shall observe the uniform
convergence on compacts of reflectionless potentials via
the convergence of $P^\sigma$'s.
The surjectivity of $\psi:\mathcal{G}\to\Xi$ will be seen
there.

\section{Reflectionless potentials}
\label{sec:2}

We start this section by reviewing the result in
\cite{ikeda-taniguchi}. 
In what follows, every element in $\mathbf{R}^n$ is regarded
as a column vector, and ${}^t\!A$ stands for the transpose of
matrix $A$.

Let 
$$
 \Sigma_0=\Bigl\{
  \sum_{j=1}^n c_j^2\delta_{p_j}\Big|
   c_j>0,\,p_j\in \mathbf{R},\,p_i\ne p_j\,(i\ne j),
   \,n=1,2,\dots
   \Bigr\},
$$
where $\delta_p$ denotes the Dirac measure concentrated at
$p$.
For $\sigma=\sum_{j=1}^n c_j^2\delta_{p_j}\in \Sigma_0$, 
we define the $n$-dimensional Ornstein-Uhlenbeck process
$\{\xi_{\sigma}(y)\}_{y\ge0}$ and the $1$-dimensional
Gaussian process $\{X_\sigma(y)\}_{y\ge0}$ by
\begin{equation}\label{eq.ou.pr}
 \begin{array}{c}
 \displaystyle
 \xi_\sigma(y)=e^{yD_{\sigma}}
     \int_0^y e^{-zD_{\sigma}} dB(z)
     ={}^t\bigl(e^{yp_j}\int_0^y e^{-zp_j}dB^j(z)
       \bigr)_{1\le j\le n},\\
 \\
 X_\sigma(y)=\langle\mathbf{c},\xi_{\sigma}(y)\rangle,
 \end{array}
\end{equation}
where 
$\{B(y)=(B^1(y),\dots,B^n(y))\}_{y\ge0}$ is an
$n$-dimensional Brownian motion on a probability space
$(\Omega,\mathcal{F},P)$, $dB(z)$ stands for the It\^o
integral with respect to $B(z)$, $D_{\sigma}$ denotes the
$n\times n$ diagonal matrix with $p_1,\dots,p_n$ as diagonal
entries, $e^A=\sum_{j=0}^\infty A^j/j!$ for $n\times n$
matrix $A$, $\mathbf{c}={}^t(c_1,\dots,c_n)$,
and $\langle\cdot,\cdot\rangle$ is the inner product in
$\mathbf{R}^n$.
It should be mentioned that the law of $X_\sigma$ does not
depend on the order of pairs $(p_j,c_j)$'s, while
$\xi_\sigma$ does. 
It is easily seen that
\begin{eqnarray*}
  \int_\Omega X_\sigma(x)X_\sigma(y) dP
 & = &\sum_{j=1}^n \frac{c_j^2}{2p_j}\{e^{p_j(x+y)}-
        e^{p_j|x-y|}\}
 \\
 & = &\int_{\mathbf{R}}
     \frac{e^{\zeta(x+y)}-e^{\zeta|x-y|}}{2\zeta}
     \sigma(d\zeta)=R_\sigma(x,y).
\end{eqnarray*}
Hence $P^\sigma$ is realized as the induced measure of
$X_\sigma$ on $\mathcal{W}$; 
$P^\sigma=P\circ X_\sigma^{-1}$.
Note that
$$
  \frac{d}{dx}\int_{\mathcal{W}} X(x)^2 dP^\sigma
  =\int_{\mathbf{R}} e^{2\zeta x} \sigma(d\zeta).
$$
Hence $\sigma=\mu$ if $P^\sigma=P^\mu$.
Thus $\Sigma_0$ is identified with $\mathcal{G}_0$.

Let $\mathcal{S}$ be the set of all sequence 
$\{\eta_j,m_j\}_{1\le j\le n}$ of length $2n$,
$n=1,2,\dots$, of positive real numbers such that
$\eta_1<\dots<\eta_n$.
The reflectionless potential $u_{\mathbf{s}}$ with
scattering data $\mathbf{s}=\{\eta_j,m_j\}_{1\le j\le
n}\in\mathcal{S}$ is by definition the function 
$$
 u_{\mathbf{s}}(x) 
 = -2 \Bigl(\frac{d}{dx}\Bigr)^2 \log \det(I+G_{\mathbf{s}}(x)),
 \quad x\in\mathbf{R},
$$
where $G_{\mathbf{s}}(x)$ is the $n\times n$ matrix given by
$$
 G_{\mathbf{s}}(x) = \biggl( 
 \frac{\sqrt{m_im_j}e^{-(\eta_i+\eta_j)x}}{\eta_i+\eta_j}
 \biggr)_{1\le i,j\le n}.
$$
Set
$$
 \Xi_0=\{u_{\mathbf{s}}|\mathbf{s}\in \mathcal{S}\}.
$$
Solving the scattering problem for the Sturm-Liouville
operator $-(d/dx)^2+u_{\mathbf{s}}$, one obtains
scattering data $\mathbf{s}\in \mathcal{S}$ from
$u_{\mathbf{s}}$ 
(\cite{kay-moses,mar-book,marchenko}).
Thus, $\Xi_0$ and $\mathcal{S}$ are identified.
It may be interesting to recall
(\cite{mar-book,miwa-jimbo-date}) that if we set 
$$
 \mathbf{s}(t)
 =\{\eta_j,m_j\exp(-2\eta_j^3 t)\}_{1\le j\le n},
$$
then the function
$v(x,t)=-u_{\mathbf{s}(t)}(x)$ solves the KdV equation 
$$
 \frac{\partial v}{\partial t}
 = \frac32 v\frac{\partial v}{\partial x}
   +\frac14 \frac{\partial^3 v}{\partial x^3}.
$$

For $\sigma\in\Sigma_0$, without loss of generality, we may
and will assume that there exist $m\le n$ and 
$1\le j(1)<\dots<j(m)\le n$ such that 
\\
(H)~~ $|p_k|\le |p_{k+1}|$, $p_{j(\ell)}>0$, 
    $p_{j(\ell)+1}=-p_{j(\ell)}$, 
    $\#\{|p_1|,\dots,|p_n|\}=n-m$, 
\\
where $1\le k\le n-1$ and $1\le \ell\le m$.
Then, the equation $\sum_{j=1}^n c_j^2/(r-p_j^2)=1$ admits
$n-m$ roots $0<r_1<\dots<r_{n-m}$.
Define the mapping $\overline{\psi}:\Sigma_0\to \mathcal{S}$ so
that 
$\overline{\psi}(\sigma)=\{\eta_j,m_j\}_{1\le j\le n}
 \in\mathcal{S}$ is given by 
\begin{eqnarray}
 && \{\eta_1<\dots<\eta_n\}
    = \{p_{j(1)},\dots,p_{j(m)},
        r_1^{1/2},\dots,r_{n-m}^{1/2}\},
 \nonumber 
 \\
 && m_i = \left\{
    \begin{array}{ll}
      \displaystyle
      2\eta_{j(\ell)}
        \frac{c_{j(\ell)+1}^2}{c_{j(\ell)}^2}
          \prod_{k\ne j(\ell)}
            \frac{\eta_k+\eta_{j(\ell)}}{
                         \eta_k-\eta_{j(\ell)}}
          \prod_{k\ne j(\ell),j(\ell)+1} 
               \frac{p_k+\eta_{j(\ell)}}{
                          p_k-\eta_{j(\ell)}},
        \quad
      & \mbox{if }i=j(\ell),
      \\
      \displaystyle
      -2\eta_i 
          \prod_{k\ne i}
            \frac{\eta_k+\eta_i}{\eta_k-\eta_i}
          \prod_{k=1}^n
               \frac{p_k+\eta_i}{p_k-\eta_i},
       & \mbox{otherwise.}
      \end{array}\right.
   \label{eq.m_i}
\end{eqnarray}
It was seen in \cite{ikeda-taniguchi} that
\begin{eqnarray}
  \log\int_{\mathcal{W}}
    \exp\biggl(-\frac12\int_0^x X(y)^2 dy\biggr)
    dP^\sigma
  =-\frac12\log\det(I+G_{\overline{\psi}(\sigma)}(x))+
 \quad
 \nonumber
 &&
 \\
  +\frac12\log\det(I+G_{\overline{\psi}(\sigma)}(0))
  -\frac{x}{2}\sum_{j=1}^n (p_j+\eta_j), \quad x\ge0.
 &&
  \label{eq.it}
\end{eqnarray}
In particular, $\psi(P^\sigma)$ in (\ref{eq.i&t}) satisfies
that
\begin{equation}\label{eq.psi}
 \psi(P^\sigma)=u_{\overline{\psi}(\sigma)}\quad
 \mbox{on }[0,\infty)\quad\mbox{for any }
 P^\sigma\in \mathcal{G}_0. 
\end{equation}
If $u,v\in\Xi_0$ coincide on $[0,\infty)$, then so
on $\mathbf{R}$, since they are real analytic.
Thus, we may and will think of $\psi(P^\sigma)$, 
$P^\sigma\in \mathcal{G}_0$, as functions on $\mathbf{R}$,
and hence $\psi$ as a mapping of $\mathcal{G}_0$ to
$\Xi_0$.

We are now ready to state our first main result.

\begin{thm}
\label{thm.bij}
(i) $\psi:\mathcal{G}_0\to\Xi_0$ is bijective.
\\
(ii) Let $P^\sigma\in \mathcal{G}_0$ and $u=\psi(P^\sigma)$.
Represent as $\sigma=\sum_{j=1}^n c_j^2\delta_{p_j}$ and
define 
$\tilde{\sigma}=\sum_{j=1}^n c_j^2\delta_{-p_j}$.
Then it holds that
$$
 u(x)=\psi(P^{\tilde{\sigma}})(-x),\quad x\le0.
$$
\end{thm}
Due to this theorem,  $\mathcal{G}_0$ and $\Xi_0$ can be 
identified.
The theorem immediately implies that

\begin{cor}\label{cor.neg}
Let
$\tilde{\mathcal{G}}_0
  =\{Q^\sigma=(P^\sigma,P^{\tilde{\sigma}})|\sigma\in\Sigma_0\}$,
where $\tilde{\sigma}$ is defined as in
Theorem~\ref{thm.bij}.
Then the mapping $\tilde{\psi}$ defined by 
$$
 \tilde{\psi}(Q^\sigma)(x)
  =\left\{
       \begin{array}{ll}
        \psi(P^\sigma)(x),  & \mbox{if }x\ge0, 
       \\
        \psi(P^{\tilde{\sigma}})(-x),\quad & \mbox{if }x<0,
      \end{array}\right.
$$
is a bijection from $\tilde{\mathcal{G}}_0$ to $\Xi_0$.
\end{cor}

Furthermore, we have that

\begin{cor}\label{cor.brown}
Let $P^\sigma\in \mathcal{G}_0$ and $u=\psi(P^\sigma)$.
Extend the Brownian motion $\{B(y)\}_{y\ge0}$ used in
(\ref{eq.ou.pr}) to $y\le0$ so that $B(y)=B(-y)$, and define
$\xi_\sigma(y)$ and $X_\sigma(y)$ by (\ref{eq.ou.pr}) for
$y\le0$:
$$
 \xi_\sigma(y)
    =e^{yD_\sigma}\int_0^y e^{-zD_\sigma} dB(z)
    =-e^{yD_\sigma}\int_y^0 e^{-zD_\sigma} dB(z),
 \quad
 X_\sigma(y)=\langle\mathbf{c},\xi_\sigma(y)\rangle.
$$
Then it holds that
$$
 u(x)=4\Bigl(\frac{d}{dx}\Bigr)^2\log\biggl(
    \int_\Omega \exp\biggl(-\frac12
    \int_{\min\{0,x\}}^{\max\{0,x\}} X_\sigma(y)^2 dy\biggr)
    dP\biggr),\quad x\in \mathbf{R}. 
$$
\end{cor}

\begin{proof}
Let $\sigma=\sum_{j=1}^n c_j^2\delta_{p_j}$ and 
$\tilde{\sigma}=\sum_{j=1}^n c_j^2\delta_{-p_j}$.
Since $D_{\tilde{\sigma}}=-D_\sigma$, it is easily seen that
$$
 \xi_\sigma(y)=\xi_{\tilde{\sigma}}(-y), \quad y\le0.
$$
Hence $X_\sigma(y)=X_{\tilde{\sigma}}(-y)$, $y\le0$, and 
$$
 \int_x^0 X_\sigma(y)^2 dy 
 =\int_0^{-x} X_{\tilde\sigma}(y)^2 dy, \quad x\le0.
$$
Since $P^{\tilde{\sigma}}=P\circ X_{\tilde{\sigma}}^{-1}$,
in conjunction with Theorem~\ref{thm.bij}(ii), this
yields that 
$$
 u(x)=\psi(P^{\tilde{\sigma}})(-x)
 =4\Bigl(\frac{d}{dx}\Bigr)^2\log\biggl(
    \int_\Omega \exp\biggl(-\frac12 \int_x^0 X_\sigma(y)^2
    dy\biggr)dP\biggr)
$$
for $x\le0$, 
which completes the proof.
\end{proof}

\begin{proof}[Proof of Theorem~\ref{thm.bij}]
(i) 
Let 
$\mathbf{s}=\{\kappa_j,q_j\}_{1\le j\le n}\in \mathcal{S}$.
For $\lambda\in\mathbf{C}$ with $\Im\lambda\ge0$, denote by
$e^+(x;\lambda)$ and $e^-(x;-\lambda)$ the right and left
Jost solutions of 
$$
 \{-(d/dx)^2+u_\mathbf{s}\}\phi=\lambda^2\phi,
$$
respectively, i.e. $e^+(x;\lambda)$ and $e^-(x;-\lambda)$
satisfy the above ordinary differential equation and 
$e^\pm(x;\pm\lambda)\sim e^{\pm\sqrt{-1}\lambda x}$ as
$x\to\pm\infty$, where and in the sequel the symbol $\pm$
takes the same sign $+$ or $-$ simultaneously.
It was shown in \cite{lundina,marchenko} that there exist 
$\lambda_j\in C^\infty(\mathbf{R};\mathbf{R})$, 
$1\le j\le n$, such that 
$\lambda_i(x)\ne\lambda_j(x)$ if $i\ne j$ for each
$x\in\mathbf{R}$, and 
\begin{equation}\label{eq.lambda_j}
 e^\pm(x;\pm\lambda)
 =e^{\pm\sqrt{-1}\lambda x} \prod_{j=1}^n 
  \frac{\lambda-(\pm\sqrt{-1}\lambda_j(x))}{
                    \lambda+\sqrt{-1}\kappa_j}.
\end{equation}
Define $k(\alpha)$, $1\le \alpha\le n$, so that 
$|\lambda_{k(\alpha)}(0)|\le|\lambda_{k(\alpha+1)}(0)|$,
$1\le \alpha\le n-1$ and 
$\lambda_{k(\alpha)}(0)=-\lambda_{k(\alpha+1)}(0)>0$ if 
$|\lambda_{k(\alpha)}(0)|=|\lambda_{k(\alpha+1)}(0)|$.
Note that, in the latter condition,
$\lambda_{k(\alpha)}(0)$ and $\lambda_{k(\alpha+1)}(0)$ have 
signs opposite to the ones in \cite{marchenko}.
The following properties were seen in \cite{marchenko};
(A)~$\lambda_j^\prime(0)<0$, $1\le j\le n$,
(B)~for $1\le \alpha\le n$,  
either of the following two cases occurs;
(a)~$\kappa_{\alpha-1}<|\lambda_{k(\alpha)}(0)|
     <\kappa_\alpha$, or
(b)~$\lambda_{k(\alpha)}(0)=-\lambda_{k(\alpha+1)}(0)
     =\kappa_\alpha$,
where $\kappa_0=0$, 
(C) it holds that
\begin{equation}\label{eq.mar.m_i}
 \frac1{q_\alpha}
 = \frac{\kappa_\alpha^2-\lambda_{k(\alpha)}(0)^2}{
        2\kappa_\alpha(\kappa_\alpha+\lambda_{k(\alpha)}(0))^2}
   \prod_{s\ne \alpha} \biggl(
     \frac{\kappa_\alpha-\lambda_{k(s)}(0)^2}{
                  \kappa_\alpha^2-\kappa_s^2}\biggr)
   \biggl(\frac{\kappa_\alpha-\kappa_s}{\kappa_\alpha+\lambda_{k(s)}(0)}
     \biggr)^2
\end{equation}
if $\kappa_\alpha\ne|\lambda_j(0)|$ for any $j=1,\dots,n$,
and 
\begin{equation}\label{eq.mar.m_i.2}
 \frac1{q_\alpha}
 =\frac{\lambda_{k(\alpha)}^\prime(0)}{
      2\kappa_\alpha\lambda_{k(\alpha+1)}^\prime(0)}
  \frac{\kappa_{\alpha+1}-\kappa_\alpha}{\kappa_{\alpha+1}+\kappa_\alpha}
   \prod_{s\ne \alpha} \biggl(
     \frac{\kappa_\alpha-\lambda_{k(s)}(0)^2}{
                  \kappa_\alpha^2-\kappa_s^2}\biggr)
   \biggl(\frac{\kappa_\alpha-\kappa_s}{\kappa_\alpha+\lambda_{k(s)}(0)}
     \biggr)^2
\end{equation}
if $\kappa_\alpha=\lambda_{k(\alpha)}(0)$, and 
\begin{equation}\label{eq.kappa}
 \prod_{j=1}^n (z-\kappa_j^2)
 = \biggl\{\prod_{j=1}^n (z-\lambda_j(0)^2)\biggr\}
   \biggl\{ 1-\sum_{j=1}^n 
    \frac{-\lambda_j^\prime(0)}{z-\lambda_j(0)^2}\biggr\}.
\end{equation}

Let $u=u_{\mathbf{s}}\in\Xi_0$ with
$\mathbf{s}=\{\kappa_j,q_j\}_{1\le j\le n}\in \mathcal{S}$.
Define
$$
 p_\alpha(\mathbf{s})=\lambda_{k(\alpha)}(0),~~
 c_\alpha(\mathbf{s})=\sqrt{-\lambda_{k(\alpha)}^\prime(0)},~~
 \sigma(\mathbf{s})=\sum_{j=1}^n c_j(\mathbf{s})^2
 \delta_{p_j(\mathbf{s})}.
$$
Set 
$\overline{\psi}(\sigma(\mathbf{s}))
 =\{\eta_j,m_j\}_{1\le j\le n}$.
Since $p_j(\mathbf{s})$'s satisfy the condition (H), by 
(\ref{eq.kappa}), we see that $\eta_j=\kappa_j$, 
$1\le j\le n$. 
Substituting these into (\ref{eq.mar.m_i}) and
(\ref{eq.mar.m_i.2}), and then comparing with
(\ref{eq.m_i}), we obtain that 
$m_j=q_j$, $1\le j\le n$.
Hence $\overline{\psi}(\sigma(\mathbf{s}))=\mathbf{s}$.
Due to (\ref{eq.psi}),
$\psi(P^{\sigma(\mathbf{s})})=u_{\mathbf{s}}$, which means
that $\psi$ is surjective.

Let $\sigma=\sum_{j=0}^\infty
c_j^2\delta_{p_j}\in\Sigma_0$, and assume that (H) is
satisfied. 
Let $\mathbf{s}=\overline{\psi}(\sigma)$.
It was shown in the proof of \cite[Lemma\,1.4]{marchenko}
that $p_j(\mathbf{s})=p_j$ and $c_j(\mathbf{s})=c_j$, 
$1\le j\le n$.
Hence, if we define the mapping $\phi:\Xi_0\to \mathcal{G}_0$ by
$\phi(u_{\mathbf{s}})=P^{\sigma(\mathbf{s})}$, then by
(\ref{eq.psi}), $\phi(\psi(P^\sigma))=P^\sigma$. 
Thus $\psi$ is injective.
\\
(ii)
Let $\sigma=\sum_{j=1}^n c_j^2\delta_{p_j}$ and
$u=\psi(P^\sigma)$.
If we set
$\overline{\psi}(\sigma)=\mathbf{s}=\{\kappa_j,q_j\}$, 
as was seen in the proof of (i), $u=u_{\mathbf{s}}$,
$p_j(\mathbf{s})=p_j$, and $c_j(\mathbf{s})=c_j$,
$j=1,\dots,n$.

Put $\tilde{u}(x)=u(-x)$, $x\in \mathbf{R}$.
Denote by $\tilde{e}^+(x;\lambda)$ and
$\tilde{e}^-(x;-\lambda)$ the right and left Jost solutions
associated with $\tilde{u}$, respectively.
It is straightforward to see that
$\tilde{e}^+(x;\lambda)=e^-(-x;-\lambda)$
and $\tilde{e}^-(x;-\lambda)=e^+(-x;\lambda)$, 
$e^+(x;\lambda)$ and $e^-(x;-\lambda)$ being the right and
left Jost solutions related with $u$, respectively.
This implies that
\begin{eqnarray*}
 && W[\tilde{e}^+(*;\lambda),\tilde{e}^-(*;-\lambda)]
   =W[e^+(*;\lambda),e^-(*;-\lambda)], 
 \\
 && W[\tilde{e}^-(*;-\xi),\tilde{e}^+(*;-\xi)]
    =W[e^-(*;\xi),e^+(*;\xi)]
\end{eqnarray*}
for any $\lambda\in\mathbf{C}$ with $\Im\lambda\ge0$ and
$\xi\in \mathbf{R}$, where $W[f,g]$ denotes the Wronskian of 
$f$ and $g$: $W[f,g]=f^\prime g-fg^\prime$.
Hence, by virtue of the direct and inverse scattering theory
(cf. \cite{mar-book}),
$\tilde{u}\in\Xi_0$ and there exist
$\tilde{q}_1,\dots,\tilde{q}_n>0$ so that, if we set 
$\tilde{\mathbf{s}}=\{\kappa_j,\tilde{q}_j\}$ then 
$\tilde{u}=u_{\tilde{\mathbf{s}}}$.
Due to (\ref{eq.lambda_j}), we have that
$$
 \tilde{e}^\pm(x;\pm\lambda)
 =e^{\pm\sqrt{-1}\lambda x}\prod_{j=1}^n
  \frac{\lambda-(\pm\sqrt{-1}(-\lambda_j(-x)))}{
                      \lambda+\sqrt{-1}\kappa_j}.
$$
By the definition of $p_j(\mathbf{s})$ and
$c_j(\mathbf{s})$, this implies that 
$p_j(\tilde{\mathbf{s}})=-p_j(\mathbf{s})=-p_j$ and
$c_j(\tilde{\mathbf{s}})=c_j(\mathbf{s})=c_j$,
$j=1,\dots,n$. 
In particular, $\sigma(\tilde{\mathbf{s}})=\tilde{\sigma}$.
Thus 
$\overline{\psi}(\tilde{\sigma})=\tilde{\mathbf{s}}$,
and hence $\tilde{u}=\psi(P^{\tilde{\sigma}})$ on
$[0,\infty)$, which completes the proof.
\end{proof}

\section{The Brownian sheet}
\label{sec:3}

\subsection{Wiener integral with respect to the Brownian sheet}

Let $\{W(p,x)\}_{(p,x)\in\mathbf{R}_+^2}$ be the Brownian 
sheet on a probability space $(\Omega,\mathcal{F},P)$, where 
$\mathbf{R}_+^2=[0,\infty)^2$, 
i.e. $\{W(p,x)\}_{(p,x)\in\mathbf{R}_+^2}$ 
is a centered Gaussian system with covariance function 
$\int_\Omega W(p,x)W(q,y)dP=\min\{p,q\}\min\{x,y\}$.
Denote by $L^2(\mathbf{R}_+^2)$ and $L^2(P)$ the spaces of 
square integrable functions with respect to the Lebesgue
measure on $\mathbf{R}_+^2$ and $P$, respectively.
There exists a linear isometry 
$\mathcal{I}:L^2(\mathbf{R}_+^2)\to L^2(P)$ such that
$$
 \mathcal{I}(\chi_{[a,b)\times[c,d)})=W(b,d)-W(a,d)-W(b,c)+W(a,c),
$$
for any $0\le a<b<\infty$ and $0\le c<d<\infty$, 
where $\chi_A$ is the indicator function of $A$.
In the sequel, we shall write 
$$
  \int_{\mathbf{R}_+^2} h(q,z)W(dq,dz)
$$
for $\mathcal{I}(h)$, and call it the Wiener integral of
$h$.

We shall see the dependence of the Wiener integrals on
parameters.
To do this, let $T>0$ and take a family 
$\phi=\{\phi(\cdot,\cdot;t)\,|\,t\in[0,T]\}
 \subset L^2(\mathbf{R}_+^2)$ such that 
\begin{equation}\label{eq.k_phi}
 K_\phi\equiv\sup_{0\le s<t\le T} \frac{1}{|t-s|}
 \int_{\mathbf{R}_+^2} |\phi(q,z;t)-\phi(q,z;s)|^2 dqdz
 <\infty,
\end{equation}
and put 
$Z_\phi(y)=\int_{\mathbf{R}_+^2} \phi(q,z;y)W(dq,dz)$,
$y\in[0,T]$.
It then holds that
\begin{equation}\label{t.conti.30}
 \int_{\Omega} |Z_\phi(t)-Z_\phi(s)|^{2m} dP  
 \le \frac{(2m)!}{2^m m!} K_\phi^m|t-s|^m
 \quad\mbox{for any }t,s\in[0,T],
\end{equation}
because, for any $h\in L^2(\mathbf{R}_+^2)$, its Wiener
integral is a centered Gaussian random variable with
variance $\|h\|_{L^2(\mathbf{R}_+^2)}^2$ and hence 
$$
 \int_{\Omega} \biggl( \int_{\mathbf{R}_+^2}
       h(q,z)W(dq,dz)\biggr)^{2m} dP 
 = \frac{(2m)!}{2^m m!}\|h\|_{L^2(\mathbf{R}_+^2)}^{2m}, 
 \quad m\in\mathbf{N}.
$$
By Kolmogorov's continuity theorem, 
$\{Z_\phi(y)\}_{y\in[0,T]}$ admits a continuous version, say 
$\{Z_\phi(y)\}_{y\in[0,T]}$ again.
We moreover have that

\begin{thm}\label{thm.conti}
Let $T>0$ and $m\in \mathbf{N},\ge2$.
Then, there exists a constant $C_{m,T}>0$ such that, for any
family 
$\phi=\{\phi(\cdot,\cdot;t)\,|\,t\in[0,T]\}
 \subset L^2(\mathbf{R}_+^2)$ with $K_\phi<\infty$, where
$K_\phi$ is defined by (\ref{eq.k_phi}),
the Wiener integral 
$$
 Z_\phi(y)=\int_{\mathbf{R}_+^2} \phi(q,z;y)W(dq,dz)
$$
satisfies that
\begin{equation}\label{t.conti.2}
 \int_{\Omega} \sup_{0\le s<t\le T} 
    \frac{|Z_\phi(t)-Z_\phi(s)|^{2m}}{|t-s|^{m-(3/2)}} \,dP
   \le C_{m,T} K_\phi^m.
\end{equation}
Moreover, if $Z_\phi(0)=0$ in addition, then it holds that
$$
 \int_{\Omega} \sup_{y\in[0,T]}|Z_\phi(y)|^{2m} dP
   \le C_{m,T}K_\phi^mT^{m-(3/2)}.
$$
\end{thm}

\begin{proof}
To see the assertion, we apply the following inequality,
which can be concluded easily from
\cite[Theorem~2.1.3]{stroock-varadhan};
for each $\alpha>0$, $\beta>2$, $T>0$, and continuous
function $f:[0,T]\to \mathbf{R}$, it holds that
$$
 \sup_{0\le s<t\le T}
   \frac{|f(t)-f(s)|^\alpha}{|t-s|^{\beta-2}}
 \le 2^{3 \alpha+2}
     \biggl(\frac{\beta}{\beta-2}\biggr)^{\alpha}
    \int_0^T\int_0^T
     \frac{|f(t)-f(s)|^\alpha}{|t-s|^\beta} dt ds.
$$
Plugging (\ref{t.conti.30}) into this estimation
with $\alpha=2m$ and $\beta=m+(1/2)$, we have that
\begin{eqnarray*}
 && \int_{\Omega} \sup_{0\le s<t\le T}
   \frac{|Z_\phi(t)-Z_\phi(s)|^{2m}}{|t-s|^{m-(3/2)}}
   dP 
 \\
 && \qquad
   \le 2^{6m+2}\biggl(\frac{2m+1}{2m-3}\biggr)^{2m}
     \frac{(2m)!}{2^m m!} K_\phi^m
     \int_0^T\int_0^T |t-s|^{-1/2}dtds.
\end{eqnarray*}
Thus we obtain (\ref{t.conti.2}).
The last inequality is an immediate consequence
of (\ref{t.conti.2}).
\end{proof}

\subsection{Representation with the Brownian sheet}

We first reconstruct $P^\sigma\in \mathcal{G}_0$ by using the
Brownian sheet.
For this purpose, let $\mathcal{Q}$ be the set of all
sequence $\alpha=\{(p_j,d_j)\}_{1\le j\le n}$ of points in
$\mathbf{R}^2$ with $p_i\ne p_j$ if $i\ne j$, $n=1,2,\dots$
Every $\sigma=\sum_{j=1}^n c_j^2\delta_{p_j}\in\Sigma_0$
determines the element 
$\{(p_j,c_j)\}_{1\le j\le n}\in \mathcal{Q}$, denoted by
$\sigma$ again, if we order $p_j$'s so that the condition
(H) is fulfilled.

For $\alpha=\{(p_j,d_j)\}_{1\le j\le n}\in \mathcal{Q}$,
$a\ge0$ and $b\in \mathbf{R}$ with $-a\le b<p_1$, 
define $0\le q_0<q_1<\dots<q_n$ by 
\begin{equation}\label{eq.def.q_n}
 q_0=b+a,\quad
 q_k=q_0+\sum_{j=1}^k |p_j-p_{j-1}|,
 \quad k=1,\dots,n \quad(p_0=b).
\end{equation}
The $\mathbf{R}^n$-valued process 
$$
 W_\alpha(y)
 = \biggl( \frac{W(q_j,y)-W(q_{j-1},y)}{
            \sqrt{q_j-q_{j-1}}} \biggr)_{1\le j\le n}
$$
is an $n$-dimensional Brownian motion, and then using this
for $\{B(z)\}$ in (\ref{eq.ou.pr}), we define 
$$
 \xi_{a,b,\alpha}(y)
 = e^{yD_\alpha}\int_0^y e^{-zD_\alpha} 
   dW_\alpha(z)
 \quad\mbox{and}\quad
 X_{a,b,\alpha}(y)
 =\langle \mathbf{d},
     \xi_{a,b,\alpha}(y)\rangle,
$$
where $D_\alpha$ denotes the diagonal matrix with $p_j$'s as
diagonal elements and $\mathbf{d}={}^t\!(d_1,\dots,d_n)$.
Then it is easily seen that
\begin{equation}\label{eq.Xabq}
  X_{a,b,\alpha}(y)
  =\int_{\mathbf{R}_+^2} h_{a,b,\alpha}(q,z;y)W(dq,dz),
\end{equation}
where
$$
 h_{a,b,\alpha}(q,z;y)
 =\sum_{j=1}^n 
  \frac{e^{(y-z)p_j}d_j}{\sqrt{q_j-q_{j-1}}}
  \chi_{[q_{j-1},q_j)\times[0,y)}(q,z).
$$
Moreover, if $\sigma\in\Sigma_0$, then, by virtue of the 
observation made in Sect.~\ref{sec:2}, it holds that
$$
  P^{\sigma}=P\circ X_{a,b,\sigma}^{-1}.
$$

We next introduce another compound Ornstein-Uhlenbeck process.
For $a\ge0$ and a piecewise continuous function 
$g:[0,\infty)\to\mathbf{R}$ with compact support, we define
$h_{a,g}(\cdot,\cdot;y)\in L^2(\mathbf{R}_+^2)$,
$y\in[0,\infty)$, by 
$$
 h_{a,g}(q,z;y)
 = e^{(y-z)(q-a)}g(q) \chi_{[0,y)}(z),
 \quad (q,z)\in\mathbf{R}_+^2,
$$
and then put
$$
 X_{a,g}(y)=\int_{\mathbf{R}_+^2}
   h_{a,g}(q,z;y) W(dq,dz),
 \quad y\in[0,\infty).
$$

We shall give some remarks on $X_{a,b,\alpha}$ and
$X_{a,g}$.
Firstly notice that $X_{a,b,\alpha}$ and $X_{a,g}$ are both
continuous Gaussian processes starting at $0$ at time $0$.
Namely, being Gaussian processes follows from their
definition by Wiener integrals.
The continuity is a consequence of the observation made
before Theorem~\ref{thm.conti} and the next lemma.

\begin{lem}\label{lem.comp.ou.est}
Let $g,\alpha$ be as above and $T>0$.
Set 
\begin{eqnarray*}
 && K_{\alpha,T} 
 = e^{2TM(\alpha)}
      \{1+T^2 M(\alpha)^2 \}S(\alpha),
 \\
 && K_{a,g,T} = \{1+(T_0+a)^2T^2\} e^{2T(T_0+a)}
    \int_0^\infty g(q)^2 dq,
\end{eqnarray*}
where 
$M(\alpha)=\sup_{1\le j\le n}|p_j|$, 
$S(\alpha)=\sum_{j=1}^n d_j^2$, and
$T_0$ is chosen so that $g(q)=0$ if $q\ge T_0$.
Then it holds that
$$ 
 K_{h_{a,b,\alpha}}\le K_{\alpha,T}
 \quad\mbox{and}\quad
 K_{h_{a,g}}\le K_{a,g,T},
$$
where $K_{h_{a,b,\alpha}}$ and $K_{h_{a,g}}$ are defined by 
(\ref{eq.k_phi}) with $\phi=h_{a,b,\alpha}$ and $h_{a,g}$,
respectively.
\end{lem}

\begin{proof}
For any $0\le s<t\le T$, it holds that
\begin{eqnarray*}
 &&
 |h_{a,b,\alpha}(q,z;t)-h_{a,b,\alpha}(q,z;s)|
  \le
  \sum_{j=1}^n
   \frac{e^{TM(\alpha)}|d_j|}{\sqrt{q_j-q_{j-1}}} 
   \chi_{[q_{j-1},q_j)\times[s,t)}(q,z)
 \\
 && \hspace{100pt}
   + \sum_{j=1}^\infty 
     \frac{M(\alpha)e^{TM(\alpha)}(t-s)|d_j|}{
                             \sqrt{q_j-q_{j-1}}}
     \chi_{[q_{j-1},q_j)\times[0,s)}(q,z),
 \\
 &&
 |h_{a,g}(q,z;t)-h_{a,g}(q,z;s)|
 \\
 && \qquad
 \le e^{T(T_0+a)}|g(q)|
     \{\chi_{[s,t)}(z)+(t-s)(T_0+a)\chi_{[0,s)}(z)\}.
\end{eqnarray*}
These imply the desired conclusion.
\end{proof}

Secondly, observe that for $\sigma,\mu\in\Sigma_0$, if $A$
and $B$ are chosen so that $A+B$ is sufficiently large,
then
\begin{equation}\label{eq.sum}
 P^{\sigma+\mu}=P\circ\{X_{a,b,\sigma}+X_{A,B,\mu}\}^{-1}.
\end{equation}
Namely, note that $h_{A,B,\mu}(q,z;y)=0$ if $q\le A+B$. 
Hence, if $A+B$ is so large that $q_n\le A+B$, where $q_n$
is defined by (\ref{eq.def.q_n}) for $\sigma$, then 
$h_{a,b,\sigma}h_{A,B,\mu}=0$, and which implies the
independence of $X_{a,b,\sigma}$ and $X_{A,B,\mu}$.
Then 
\begin{eqnarray*}
 && \int_\Omega \{X_{a,b,\sigma}(x)+X_{A,B,\mu}(x)\}
    \{X_{a,b,\sigma}(y)+X_{A,B,\mu}(y)\} dP
 \\
 && \qquad
    = R_\sigma(x,y)+R_\mu(x,y)
    = R_{\sigma+\mu}(x,y).
\end{eqnarray*}
Thus $P^{\sigma+\mu}$ is realized as the law of 
$X_{a,b,\sigma}+X_{A,B,\mu}$.

Thirdly, if $\sigma\in\Sigma$ is of the form 
$$
 \sigma(d\xi)=f(\xi)d\xi+\mu(d\xi),
$$
where $f:\mathbf{R}\to[0,\infty)$ is a piecewise continuous
function with compact support and $\mu\in\Sigma_0$, then,
choosing $a>0$ so that $\mbox{\rm supp}\,f\subset[-a,a]$,
and setting $g(\xi)=\sqrt{f(\xi-a)}$, we have that
\begin{equation}\label{eq.g+mu}
 P^\sigma=P\circ\{X_{a,g}+X_{A,B,\mu}\}^{-1}
\end{equation}
for $A$ and $B$ with sufficiently large $A+B$.
In fact, it holds that
$$
 \sigma(d\xi)=g(\xi+a)^2\chi_{[-a,\infty)}(\xi)d\xi
       +\mu(d\xi),
$$
and we may and will think of $g$ as a piecewise continuous
function on $[0,\infty)$ with compact support.
It is easily seen that the covariance function of $X_{a,g}$
is 
$$
  \int_{\Omega} X_{a,g}(x)X_{a,g}(y)dP
  =\int_{\mathbf{R}} \frac{e^{\xi(x+y)}-e^{\xi|x-y|}}{2\xi}
   g(\xi+a)^2\chi_{[-a,\infty)}(\xi) d\xi.
$$
Take $\gamma>0$ so that 
$\mbox{\rm supp}\,\mu\subset[-\gamma,\gamma]$.
Since $\mbox{\rm supp}\,g\subset[0,2a]$, 
for $A\ge0$ and $B\le0$ such that 
$-A\le B\le-\gamma$ and $2a<A+B$, we have that
$h_{a,g}h_{A,B,\mu}=0$. 
Then $X_{a,g}$ and $X_{A,B,\mu}$ are independent, and hence
the Gaussian process $X_{a,g}+X_{A,B,\mu}$ possesses the
covariance function $R_\sigma(x,y)$.
Thus $P^\sigma$ coincides with the law of
$X_{a,g}+X_{A,B,\mu}$.

Finally, Theorem~\ref{thm.conti} and
Lemma~\ref{lem.comp.ou.est} yields that

\begin{prop}\label{prop.comp.ou}
Let $g$, $\alpha$, $a,b$ be as above.
Then, for any $T>0$ and $m\in \mathbf{N}$, there exists a
constant $C_{m,T}$, depending only on $T$ and $m$, such that 
the following estimations hold with
$(Z,K)=(X_{a,b,\alpha},K_{\alpha,T})$ or
$(Z,K)=(X_{a,g},K_{a,g,T})$.

\begin{eqnarray*}
 && \int_{\Omega} \sup_{0\le s<t\le T} 
    \frac{|Z(t)-Z(s)|^{2m}}{|t-s|^{m-(3/2)}} \,dP
   \le C_{m,T} K^m,
\\
 && \int_{\Omega} \sup_{y\in[0,T]}|Z(y)|^{2m} dP
   \le C_{m,T}K^mT^{m-(3/2)}.
\end{eqnarray*}
\end{prop}

\section{Generalized reflectionless potentials}
\label{sec:4}

In this section, we shall show that the convergence of
$P^\sigma\in \mathcal{G}_0$ implies that of reflectionless
potentials to generalized one in the topology of uniform
convergence on compacts.

For $T>0$, let $\mathcal{W}_T$ be the space of all
continuous $w:[0,T]\to \mathbf{R}$ with $w(0)=0$.
Naturally $\mathcal{W}_T\subset \mathcal{W}$, and every
probability measure $P$ can be restricted to
$\mathcal{W}_T$.
The restriction will be denoted by $P|_{\mathcal{W}_T}$.
For $\sigma\in\Sigma$, put
$$
 \Phi_\sigma(x)=\int_{\mathcal{W}} \exp\biggl(
  -\frac12\int_0^x X(y)^2 dy\biggr)dP^\sigma.
$$
As will be seen in the next theorem, $\Phi_\sigma$ is $C^2$,
and then one can define
$$
 \psi(P^\sigma)=4\Bigl(\frac{d}{dx}\Bigr)^2\log\Phi_\sigma.
$$

Our goal of this section is

\begin{thm}\label{thm.conv.gen}
(i) 
For  $\sigma\in\Sigma$, $\Phi_\sigma$ is $C^2$.
\\
(ii)
Let $\sigma_n\in\Sigma_0$ and $\sigma\in\Sigma$.
Suppose that
$\bigcup_{n\in \mathbf{N}}\mbox{\rm supp}\,\sigma_n
  \subset[-\beta,\beta]$ for some $\beta>0$, and $\sigma_n$
tends to $\sigma$ vaguely.
Then $\Phi_{\sigma_n}$ and its first and second derivatives
$\Phi_{\sigma_n}^\prime$ and 
$\Phi_{\sigma_n}^{\prime\prime}$ converge to 
$\Phi_\sigma$, $\Phi_\sigma^\prime$, and 
$\Phi_\sigma^{\prime\prime}$ uniformly on every bounded
interval in $[0,\infty)$, respectively.
In particular, $\psi(P^{\sigma_n})$ tends to
$\psi(P^\sigma)$ uniformly on every bounded
interval in $[0,\infty)$.
Moreover, for every $\varepsilon>0$, there exists
$n_0\in\mathbf{N}$ such that 
\begin{equation}\label{t.conv.gen.1}
  \mbox{\rm Spec}(-(d/dx)^2+\psi(P^{\sigma_n}))
 \subset[-\beta^2-\sigma(\mathbf{R})-\varepsilon,\infty),
 \quad n\ge n_0.
\end{equation}
Finally, there exists $u\in\Xi$ such that
$\psi(P^\sigma)=u$ on $[0,\infty)$.
\\
(iii)
Let $g_n:\mathbf{R}\to[0,\infty)$ be piecewise continuous,
and $\mu\in\Sigma_0$.
Assume that 
$$
 \bigcup_{n\in \mathbf{N}}
 \mbox{\rm supp}\,g_n\subset[-\beta,\beta]
 \quad\mbox{for some }\beta>0, \quad
 \sup_{n\in\mathbf{N}}\int_{\mathbf{R}}g_n(\xi)^2d\xi<\infty,
$$
and $\sigma_n\in\Sigma$ defined by
$\sigma_n(d\xi)=g_n(\xi)^2d\xi+\mu(d\xi)$ converges to some
$\sigma\in\Sigma$ vaguely.
Then $\Phi_{\sigma_n}$,
$\Phi_{\sigma_n}^\prime$, and 
$\Phi_{\sigma_n}^{\prime\prime}$ converge to 
$\Phi_\sigma$, $\Phi_\sigma^\prime$, and 
$\Phi_\sigma^{\prime\prime}$ uniformly on every bounded
interval in $[0,\infty)$, respectively.
In particular, $\psi(P^{\sigma_n})$ tends to
$\psi(P^\sigma)$ uniformly on every bounded
interval in $[0,\infty)$.
\\
(iv) For every $u\in\Xi$, there exists $P^\sigma\in
\mathcal{G}$ such that $\psi(P^\sigma)=u$ on $[0,\infty)$.
\end{thm}

We shall give several remarks on the theorem before getting
into the proof.

\begin{rem}\label{rem.gen.rl}
(a) Repeating the arguments in Lemmas~\ref{lem.derivative},
\ref{lem.conv}, and \ref{lem.i&iii} below, one can show that
$\Phi_\sigma$ is $C^\infty$.

\noindent
(b) 
Let $\sigma\in\Sigma$.
Fix $\beta>0$ so that 
$\mbox{\rm supp}\,\sigma\subset[-\beta,\beta]$, and define
$\sigma_n\in\Sigma_0$ by 
$\sigma_n(d\xi)=\sum_{j=-n}^n
   \sigma([j\beta/n,(j+1)\beta/n)) \delta_{j\beta/n}$.
Then $\sigma_n$'s satisfy the assumption in (ii).

\noindent
(c) The identification of $\Xi_0$ and $\mathcal{G}_0$
extends to that of $\Xi$ and $\mathcal{G}$ as follows.
First let $P^\sigma\in \mathcal{G}$.
Define $\sigma_n\in\Sigma_0$ as in (b).
By Theorem~\ref{thm.bij},
$\psi(P^{\sigma_n})(x)=\psi(P^{\tilde{\sigma}_n})(-x)$,
$x\le0$. 
Define $\tilde{\sigma}\in\Sigma$ by 
$\tilde{\sigma}(A)=\sigma(-A)$, 
$A\in \mathcal{B}(\mathbf{R})$, where $-A=\{-x|x\in A\}$.
Since $\mbox{\rm supp}\,\tilde{\sigma}_n\subset
 [-\beta,\beta]$ and $\tilde{\sigma}_n$ tends to
$\tilde\sigma$ vaguely, by (ii), we see that
$\psi(P^{\tilde{\sigma}_n})$ converges to 
$\psi(P^{\tilde{\sigma}})$ uniformly on compacts in
$[0,\infty)$.
As will be seen in the proof of Lemma~\ref{lem.conv} below,
there exist $u\in\Xi$ and a subsequence $\{\sigma_{n_j}\}$
of $\{\sigma_n\}$ such that $\psi(P^{\sigma_{n_j}})$
converges to $u\in\Xi$ uniformly on compacts in
$\mathbf{R}$. 
Hence we have that
$u=\psi(P^{\sigma})$ on $[0,\infty)$ and
$=\psi(P^{\tilde{\sigma}})$ on $(-\infty,0]$.

Conversely, let $u\in\Xi$.
As will be seen in the proof of (iv) (Lemma~\ref{lem.v}
below), there exist $P^{\sigma_n}\in \mathcal{G}_0$,
$n\in\mathbf{N}$, such that
$\psi(P^{\sigma_n})$ converges to $u$ uniformly on compacts
in $\mathbf{R}$, 
$\bigcup_{n\in\mathbf{N}} \mbox{\rm supp}\,\sigma_n
 \subset[-\beta,\beta]$ for some $\beta>0$,
and $\sigma_n$ tends to some
$\sigma\in\Sigma$ vaguely.
Then, in repetition of the above argument, we see that
$u$ coincides with $\psi(P^{\sigma})$ on $[0,\infty)$ and
$\psi(P^{\tilde{\sigma}})$ on $(-\infty,0]$.

\noindent
(d) A correspondence between $\Xi$ and $\Sigma$ was studied
by Marchenko \cite{marchenko} and Kotani \cite{kotani} in
an analytical manner. 
The relation between $\Xi$ and $\mathcal{G}$ investigated
above is a probabilistic counterpart to their observation.

\noindent
(e) 
Every $\psi(P^\sigma)$, $P^\sigma\in \mathcal{G}$, can be
approximated by $\psi(P^{\sigma_n})$'s with $\sigma_n$ of
the form as described in (iii).
Namely, let $P^\sigma\in \mathcal{G}$.
Take a nonnegative $C^\infty$ function 
$\phi:\mathbf{R}\to \mathbf{R}$ with compact support such
that $\int_{\mathbf{R}} \phi(x)dx=1$.
Define $g_n:\mathbf{R}\to[0,\infty)$ by
$g_n(x)^2=\int_{\mathbf{R}} n\phi(n(x-\xi))\sigma(d\xi)$,
and set $\sigma_n(d\xi)=g_n(\xi)^2d\xi$.
Then $\bigcup_{n\in\mathbf{N}}\,g_n\subset[-\beta,\beta]$
for some $\beta>0$, 
$\int_{\mathbf{R}} g_n(\xi)^2 d\xi=\sigma(\mathbf{R})$, 
and $\sigma_n$ converges to $\sigma$ vaguely.

\noindent
(f) The convergence discussed in (iii) relates to the
convergence of finite-zone potentials to reflectionless
ones discussed in \cite{dmn,novikov}.
Namely, for $u\in\Xi$ of finite-zone, the $\sigma$ appearing 
in (iv) was computed by Kotani \cite{kotani} to be
represented as $\sigma(d\xi)=g(\xi)^2d\xi+\mu(d\xi)$ for
some piecewise continuous $g$ with compact support and 
$\mu\in\Sigma_0$.
As finite-zone potentials tends to a reflectionless
potential, the support of $g$ shrinks to a discrete point
set (\cite{novikov}).
This is the situation investigated in (iii).
\end{rem}

We now proceed to the proof of Theorem~\ref{thm.conv.gen}.
It is broken into several steps, each step being a lemma. 
In the sequel, let $\{W(p,x)\}_{(p,x)\in \mathbf{R}_+^2}$ be
the Brownian sheet on $(\Omega,\mathcal{F},P)$ as in
Sect.~\ref{sec:3}.

\begin{lem}
\label{lem.C1}
Let $T>0$ and
$\{\{Z_\beta(y)\}_{y\in[0,T]}\,|\,\beta\in\Lambda\}$ be a 
family of continuous processes $Z_\beta$ defined on
$(\Omega,\mathcal{F},P)$ with $Z_\beta(0)=0$.
Suppose that
$$
 A_m=\sup_{\beta\in\Lambda} \int_{\Omega} 
  \sup_{0\le s<t\le T} 
   \frac{|Z_\beta(t)-Z_\beta(s)|^{2m}}{|t-s|^{m-(3/2)}}
 dP<\infty,
 \quad m=2,3,\dots
$$
Then the family $\{P\circ
Z_\beta^{-1}\}_{\beta\in\Lambda}$ of the laws of $Z_\beta$'s
on $\mathcal{W}_T$ is tight.

Let $Q:\mathbf{R}^2\to \mathbf{R}$ be a polynomial, and put
$$
 \Phi_{\beta,\gamma}(x)=\int_{\Omega}
  Q(Z_\beta(x),Z_\gamma(x))\exp\biggl(-\frac12 
  \int_0^x Z_\beta(y)^2 dy \biggr) dP,
 \quad x\in[0,T].
$$
Then $\Phi_{\beta,\gamma}$, $\beta,\gamma\in\Lambda$, 
are equi-continuous and uniformly bounded on $[0,T]$.

Finally, if $Q\equiv1$, then $\Phi_{\beta,\gamma}$'s are all
$C^1$, and $\Phi_{\beta,\gamma}^\prime$'s are also
equi-continuous and uniformly bounded on $[0,T]$.
\end{lem}

\begin{proof}
The finiteness of $A_m$ implies the tightness.
It also yields that
$$
 B_m=\sup_{\beta\in\Lambda}\int_{\Omega}
      \sup_{t\in[0,T]} |Z_\beta(t)|^{2m} dP<\infty, \quad
 m=2,3,\dots
$$
Then, as an application of the dominated convergence theorem
and the second assertion, we obtain the third assertion. 

To see the second assertion, let $k$ be the degree of $Q$
and take $C_0<\infty$ such that 
\begin{eqnarray*}
 & |Q(a,b)-Q(c,d)|\le C_0(1+|a|+|b|+|c|+|d|)^{k-1}
   (|a-c|+|b-d|),
 &
 \\
 & |Q(a,b)|\le C_0(1+|a|+|b|)^k,
   \quad a,b,c,d\in \mathbf{R}.&
\end{eqnarray*}
Since $|e^{-\xi}-e^{-\eta}|\le |\xi-\eta|$ for
$\xi,\eta\ge0$, we have that
\begin{eqnarray*}
 && |\Phi_{\beta,\gamma}(x)-\Phi_{\beta,\gamma}(x^\prime)|
 \\
 && \quad
   \le C_0 \biggl(\int_\Omega
      \bigl\{ 1+2\sup_{y\in[0,T]}|Z_\beta(y)|
            +2\sup_{y\in[0,T]}|Z_\gamma(y)|\bigr\}^{4(k-1)/3}
         dP\biggr)^{3/4}
 \\
 && \qquad \qquad 
    \times \biggl(\int_\omega
         \{|Z_\beta(x)-Z_\beta(x^\prime)|
           +|Z_\gamma(x)-Z_\gamma(x^\prime)|\}^4 dP 
        \biggr)^{1/4}
 \\
 && \qquad
     + \frac{C_0}2 |x-x^\prime| \int_\Omega 
         \bigl\{ 1+\sup_{y\in[0,T]}|Z_\beta(y)|
             +\sup_{y\in[0,T]}|Z_\gamma(y)|\bigr\}^{k+2}
         dP.
\end{eqnarray*}
Hence there exists a constant $C<\infty$, depending only on
$A_m$'s and $B_m$'s, such that
$$
 \sup_{\beta,\gamma\in\Lambda}
 |\Phi_{\beta,\gamma}(x)-\Phi_{\beta,\gamma}(x^\prime)|
 \le C\{|x-x^\prime|^{1/8}+|x-x^\prime|\},
 \quad x,x^\prime\in[0,T].
$$
Thus $\Phi_{\beta,\gamma}$'s are equi-continuous on $[0,T]$.
Since $\Phi_{\beta,\gamma}(0)=Q(0,0)$,
$\Phi_{\beta,\gamma}$'s are then uniformly bounded on
$[0,T]$. 
\end{proof}

\begin{lem}\label{lem.derivative}
Let $\sigma=\sum_{j=1}^n c_j^2\delta_{p_j}\in\Sigma_0$ and
$-a\le b< p_1$.
$\Phi_\sigma$ is $C^\infty$ and its first and second
derivatives are represented as 
\begin{eqnarray*}
 && \Phi_\sigma^\prime(x)
   =-\frac12 \int_{\Omega}
      X_{a,b,\sigma}(x)^2 
     \exp\biggl( -\frac12\int_0^x 
     X_{a,b,\sigma}(y)^2 dy \biggr)
    dP,
 \\
 &&
   \Phi_\sigma^{\prime\prime}(x)
   =-\frac14 \int_{\Omega} \biggl\{ 
      2\sigma(\mathbf{R})
      +4X_{a,b,\sigma}(x)
        X_{a,b,\alpha(\sigma)}(x)
      -X_{a,b,\sigma}(x)^4\biggr\}
  \\
  && \hspace{175pt}
    \times
     \exp\biggl( -\frac12\int_0^x 
     X_{a,b,\sigma}(y)^2 dy \biggr)
    dP,
\end{eqnarray*}
where $\alpha(\sigma)=\{(p_j,p_jc_j)\}_{1\le j\le n}$.
\end{lem}

\begin{proof}
By (\ref{eq.it}), $\Phi_\sigma$ is $C^\infty$.
Since $P^\sigma=P\circ X_{a,b,\sigma}^{-1}$, 
$$
 \Phi_\sigma(x)=\int_{\Omega}\exp\biggl(-\frac12
  \int_0^x X_{a,b,\sigma}(y)^2 dy \biggr)dP.
$$
Moreover, by Proposition~\ref{prop.comp.ou}, we have that 
$$
 \int_\Omega \sup_{y\in[0,T]}
  |X_{a,b,\sigma}(y)|^{2m} dP < \infty, \quad T>0, m\ge2. 
$$
By an application of the dominated convergence theorem,
the desired expression of the first derivative is obtained.

Rewrite $\xi_{a,b,\sigma}$ used to define $X_{a,b,\sigma}$
as
$$
 \xi_{a,b,\sigma}(y)
 = W_\alpha(y)+\int_0^y 
    D_\sigma\xi_{a,b,\sigma}(z) dz.
$$
Then, as an application of It\^o's formula, we have that
\begin{eqnarray*}
 X_{a,b,\sigma}(x)^2
 &=& 2\int_0^x X_{a,b,\sigma}(z)
    \langle \mathbf{c},dW_\alpha(z)\rangle
 \\
 && \qquad
 + \int_0^x \biggl\{ \sum_{j=1}^n c_j^2 
    +2X_{a,b,\sigma}(z)
     X_{a,b,\alpha(\sigma)}(z) \biggr\} dz,
\end{eqnarray*}
and hence that
\begin{eqnarray*}
  \Phi_\sigma^\prime(x)
  =  -\frac14\int_{\Omega}\!\int_0^x \!
    \biggl\{ 2\sum_{j=1}^n c_j^2 
    +4X_{a,b,\sigma}(z)
      X_{a,b,\alpha(\sigma)}(z)-X_{a,b,\sigma}(z)^4
     \biggr\} 
 &&
 \\
   \times 
   \exp\biggl( -\frac12\int_0^z
   X_{a,b,\sigma}(y)^2 dy\biggr) dz
   dP.
 &&
\end{eqnarray*}
This implies that
$\Phi_{a,b,\sigma}^\prime$ is
continuously differentiable and the second derivative of 
$\Phi_{a,b,\sigma}$ has the desired representation, because
$\sigma(\mathbf{R})=\sum_{j=1}^n c_j^2$.
\end{proof}

\begin{lem}\label{lem.conv}
Let $\sigma_n\in\Sigma_0$ and $\sigma\in\Sigma$.
Suppose that 
$\bigcup_{n\in \mathbf{N}} \mbox{\rm supp}\,\sigma_n
 \subset[-\beta,\beta]$ for some $\beta>0$ and that
$\sigma_n$ tends to $\sigma$ vaguely.
Then $\Phi_\sigma$ is $C^2$, and $\Phi_{\sigma_n}$,
$\Phi_{\sigma_n}^\prime$, and 
$\Phi_{\sigma_n}^{\prime\prime}$ converge to 
$\Phi_\sigma$, $\Phi_\sigma^\prime$, and 
$\Phi_\sigma^{\prime\prime}$ uniformly on every bounded
interval in $[0,\infty)$, respectively.
Moreover, the assertion (ii) in
Theorem~\ref{thm.conv.gen} holds.
\end{lem}

\begin{proof}
Let $a\ge0$ and $b\in \mathbf{R}$ satisfy that
$-a\le b<-\beta$.
Due to the assumption, it holds that
\begin{equation}\label{l.conv.20}
 \sup_{n\in \mathbf{N}} M(\sigma_n)\le \beta
 \quad\mbox{and}\quad
 \sup_{n\in \mathbf{N}}S(\sigma_n)<\infty.
\end{equation}
By Proposition~\ref{prop.comp.ou} we see that 
\begin{equation}\label{l.conv.21}
 \sup_{n\in \mathbf{N}} \int_{\Omega} \sup_{0\le s<t\le T}
  \frac{|X_{a,b,\sigma_n}(t)-X_{a,b,\sigma_n}(s)|^{2m}}{
    |t-s|^{m-(3/2)}} dP <\infty,
 \quad T>0,m\ge2.
\end{equation}
Thus $\{P^{\sigma_n}|_{\mathcal{W}_T}\}_{n\in \mathbf{N}}$
is tight for any $T>0$.  

Since $\sigma_n$ tends to $\sigma$ vaguely and 
$\mbox{\rm supp}\,\sigma_n\subset[-\beta,\beta]$, 
$n\in\mathbf{N}$, we obtain the convergence of
$R_{\sigma_n}(x,y)$ to $R_\sigma(x,y)$ for every $x,y\ge0$.
Hence every finite dimensional distribution of
$P^{\sigma_n}$ tends to that of $P^\sigma$.
In conjunction with the tightness, this implies that
$P^{\sigma_n}|_{\mathcal{W}_T}$ converges to 
$P^\sigma|_{\mathcal{W}_T}$ weakly for any $T>0$.
In particular, $\Phi_{\sigma_n}\to\Phi_\sigma$ point wise.

Since $M(\alpha(\sigma_n))=M(\sigma_n)$ and
$S(\alpha(\sigma_n))\le \beta^2S(\sigma_n)$, by
(\ref{l.conv.20}) and Proposition~\ref{prop.comp.ou}, we
have that
\begin{equation}\label{l.conv.22}
 \sup_{n\in \mathbf{N}} \int_{\Omega} \sup_{0\le s<t\le T}
 \frac{|X_{a,b,\alpha(\sigma_n)}(t)
        -X_{a,b,\alpha(\sigma_n)}(s)|^{2m}}{
    |t-s|^{m-(3/2)}} dP <\infty,
 ~T>0,\,m\ge2.
\end{equation}
Then the equi-continuity and the uniform boundedness of
$\Phi_{\sigma_n}$, $\Phi_{\sigma_n}^\prime$, and
$\Phi_{\sigma_n}^{\prime\prime}$ on any bounded interval in
$[0,\infty)$ follow from (\ref{l.conv.21}),
(\ref{l.conv.22}), Lemmas~\ref{lem.C1} and
\ref{lem.derivative}, and the fact that
$\sigma_n(\mathbf{R})\to\sigma(\mathbf{R})$ as
$n\to\infty$. 
In conjunction with the point wise convergence of
$\Phi_{\sigma_n}$ to $\Phi_\sigma$, we see that
$\Phi_\sigma$ and that $\Phi_{\sigma_n}$,
$\Phi_{\sigma_n}^\prime$, and
$\Phi_{\sigma_n}^{\prime\prime}$ tend to $\Phi_\sigma$,
$\Phi_\sigma^\prime$, and $\Phi_\sigma^{\prime\prime}$
uniformly on any bounded interval in $[0,\infty)$,
respectively.
In particular, the first assertion of (ii) holds.

We shall show the second assertions of (ii).
By \cite[Lemma~1.4]{marchenko}, it holds that
$$
 \mbox{\rm Spec}(-(d/dx)^2+\psi(P^{\sigma_n}))
 \subset[-\lambda,\infty)
$$
for some $\lambda>0$ with
$M(\sigma_n)^2<\lambda
  \le M(\sigma_n)^2+\sigma_n(\mathbf{R})$.
Since $M(\sigma_n)^2\le \beta^2$ and 
$\sigma_n(\mathbf{R})\to\sigma(\mathbf{R})$ as
$n\to\infty$, we obtain the second assertion of (ii).

To see the last assertion of (ii), let
$u_n=\psi(P^{\sigma_n})\in\Xi_0$.
By (\ref{t.conv.gen.1}), $\{u_n\}_{n\in\mathbf{N}}$ is
precompact in the topology of uniform convergence on bounded
intervals in $\mathbf{R}$ (\cite[Lemma~2.3]{marchenko}).
Hence, by (\ref{t.conv.gen.1}), there exists $u\in\Xi$ and a
subsequence $\{u_{n_j}\}_{j\in\mathbf{N}}$ such that
$u_{n_j}$ converges to $u\in\Xi$ uniformly on any compact
interval in $\mathbf{R}$.
Combined with the convergence of $\psi(P^{\sigma_n})$ to
$\psi(P^\sigma)$ on $[0,\infty)$, we see that
$\psi(P^\sigma)=u$ on $[0,\infty)$.
\end{proof}

\begin{lem}\label{lem.i&iii}
The assertions (i) in Theorem~\ref{thm.conv.gen} holds.
\end{lem}

\begin{proof}
Let $\sigma\in\Sigma$, and define
$\sigma_n\in\Sigma_0$ by 
$$
 \sigma_n(d\xi)=\sum_{j=-n}^n 
   \sigma([j\beta/n,(j+1)\beta/n))\delta_{j\beta/n},
$$
where $\beta>0$ is chosen so that 
$\mbox{\rm supp}\,\sigma\subset[-\beta,\beta]$.
Then $\bigcup_{n\in \mathbf{N}}
      \mbox{\rm supp}\,\sigma_n\subset[-\beta,\beta]$
and $\sigma_n\to\sigma$ vaguely.
By Lemma~\ref{lem.conv}, $\Phi_\sigma$ is $C^2$.
\end{proof}

\begin{lem}\label{lem.g+mu}
Let $g:[0,\infty)\to[0,\infty)$ be piecewise continuous and
$\mu\in\Sigma_0$.
Assume that 
$\mbox{\rm supp}\,g, \mbox{\rm supp}\,\mu\subset[-\beta,\beta]$
for some $\beta>0$.
Define $\sigma\in\Sigma$ by
$$
 \sigma(d\xi)=g(\xi)^2 d\xi+\mu(d\xi).
$$
For $a>\beta$, $A>0$, and $B<0$ with 
$-A\le B\le-\beta$ and $a+\beta\le A+B$, define
$X_\sigma=X_{a,g_a}+X_{A,B,\mu}$ and
$\widetilde{X}_\sigma
 =X_{a,\tilde{g}_a}+X_{A,B,\alpha(\mu)}$,
where $g_a(x)=g(x-a)$ and $\tilde{g}_a(x)=(x-a)g(x-a)$.
Then it holds that
$$
 \Phi_\sigma^{\prime\prime}(x)
 =-\frac14\int_{\Omega} \{2\sigma(\mathbf{R})
   +4X_\sigma(x)\widetilde{X}_\sigma(x)-X_\sigma^4(x)\}
   \exp\biggl(-\frac12\int_0^x X_\sigma(y)^2dy\biggr) dP.
$$
\end{lem}

\begin{proof}
Define $\sigma_n\in\Sigma_0$ by
$$
 \sigma_n(d\xi)=\sum_{j=1}^n g_a(j(a+\beta)/n)^2
    \frac{a+\beta}{n}\delta_{(j(a+\beta)/n)-a}.
$$
Then, $\mbox{\rm supp}\,(\sigma_n+\mu)
 \subset[-a,\beta]\cup\mbox{\rm supp}\,\mu$,
$(\sigma_n+\mu)(\mathbf{R})\le 
 (a+\beta)\sup|g|^2+\mu(\mathbf{R})$, and 
$\sigma_n+\mu$ tends to $\sigma$ vaguely.
By Lemma~\ref{lem.conv},
$\Phi_{\sigma_n+\mu}^{\prime\prime}$ converges to 
$\Phi_\sigma^{\prime\prime}$ uniformly on any bounded
interval in $[0,\infty)$.

Due to (\ref{eq.sum}), we have that
$P^{\sigma_n+\mu}=P\circ\{X_{a,-a,\sigma_n}+X_{A,B,\mu}\}^{-1}$.
Moreover, in repetition of the argument used in the proof of
Lemma~\ref{lem.derivative}, we see that
\begin{eqnarray}
  \Phi_{\sigma_n+\mu}^{\prime\prime}(x)
  & = & -\frac14 \int_{\Omega} \biggl\{ 
       2\{\sigma_n(\mathbf{R})+\mu(\mathbf{R})\}
 \nonumber
 \\
 && 
      +4\{X_{a,-a,\sigma_n}(x)+X_{A,B,\mu}(x)\}
        \{X_{a,-a,\alpha(\sigma_n)}(x)
             +X_{A,B,\alpha(\mu)}(x)\}
 \nonumber
 \\
 &&\qquad\qquad\qquad
      -\{X_{a,-a,\sigma_n}(x)+X_{A,B,\mu}(x)\}^4
     \biggr\}
 \nonumber
  \\
  && \qquad
    \times
     \exp\biggl(-\frac12\int_0^x 
      \{X_{a,-a,\sigma_n}(y)+X_{A,B,\mu}(y)\}^2 
     dy \biggr) dP.
 \label{l.g+mu.21}
\end{eqnarray}
Since $h_{a,-a,\sigma_n}(*;y)$ and
$h_{a,-a,\alpha(\sigma_n)}(*;y)$ tend to $h_{a,g_a}$ and
$h_{a,\tilde{g}_a}$ in $L^2(\mathbf{R}_+^2)$ for every
$y\in[0,\infty)$, respectively,  
$X_{a,-a,\sigma_n}(y)$ and
$X_{a,-a,\alpha(\sigma_n)}(y)$ converge to $X_{a,g_a}$ and
$X_{a,\tilde{g}_a}$ in $L^2(P)$ for every $y\in[0,\infty)$,
respectively.
Moreover, by Proposition~\ref{prop.comp.ou}, we have that
$$
 \sup_{n\in \mathbf{N}}\int_{\Omega} \sup_{y\in[0,T]}
    \bigl\{ |X_{a,-a,\sigma_n}(y)|^{2m} 
          +|X_{a,-a,\alpha(\sigma_n)}(y)|^{2m} \bigr\}
  dP<\infty
$$
for any $T>0$ and $m\in\mathbf{N}$.
Then, letting $n\to\infty$ in (\ref{l.g+mu.21}), we obtain
the desired representation of $\Phi_\sigma^{\prime\prime}$. 
\end{proof}

\begin{lem}\label{lem.iv}
The assertion (iii) of Theorem~\ref{thm.conv.gen} holds.
\end{lem}

\begin{proof}
Take $a>\beta$ and $A>0$, $B<0$ so that $-A\le
B<\inf(\mbox{\rm supp}\,\mu)$, $a+\beta<A+B$, and
define $X_{\sigma_n}=X_{a,g_{n,a}}+X_{A,B,\mu}$
and
$\widetilde{X}_{\sigma_n}=X_{a,\tilde{g}_{n,a}}
 +X_{A,B,\alpha(\mu)}$, where $g_{n,a}(\xi)=g_n(\xi-a)$ and 
$\tilde{g}_{n,a}(\xi)=(\xi-a)g_n(\xi-a)$.
By (\ref{eq.g+mu}), $P^{\sigma_n}=P\circ X_{\sigma_n}^{-1}$,
and
$$
 \Phi_{\sigma_n}(x)
 =\int_\Omega \exp\biggl(-\frac12 \int_0^x
   X_{\sigma_n}(y)^2 dy \biggr) dP.
$$
Since $\sup_{n\in\mathbf{N}} 
  \int_{\mathbf{R}} \{(g_n)_a(\xi)\}^2d\xi<\infty$ and 
$\sup_{n\in\mathbf{N}} 
  \int_{\mathbf{R}} \{(\tilde{g}_n)_a(\xi)\}^2d\xi<\infty$,
it follows from Proposition~\ref{prop.comp.ou} that
$$
 \sup_{n\in \mathbf{N}}\int_{\Omega}\sup_{0\le s<t\le T}
 \frac{|X_{\sigma_n}(t)-X_{\sigma_n}(s)|^{2m}
      +|\widetilde{X}_{\sigma_n}(t)
           -\widetilde{X}_{\sigma_n}(s)|^{2m}}{
     |t-s|^{m-(3/2)}} dP <\infty
$$
for any $T>0$ and $m\ge2$.
We then obtain the desired convergence in repetition of the
proof of Lemma~\ref{lem.conv}, only this time with
$X_{\sigma_n}$, $\widetilde{X}_{\sigma_n}$, and
Lemma~\ref{lem.g+mu} for $X_{a,b,\sigma_n}$,
$X_{a,b,\alpha(\sigma_n)}$, and Lemma~\ref{lem.derivative}.
We omit the details.
\end{proof}

\begin{lem}\label{lem.v}
The assertion (iv) of Theorem~\ref{thm.conv.gen} holds.
\end{lem}

\begin{proof}
Let $u\in\Xi$ and suppose that 
$\{u_n\}_{n\in \mathbf{N}}\subset\Xi_0$ satisfies that
$u_n$ converges to $u$ uniformly on any bounded interval in
$\mathbf{R}$ and 
$\bigcup_{n\in \mathbf{N}}\mbox{\rm Spec}(-(d/dx)^2+u_n)
  \subset[-\lambda,\infty)$ for some $\lambda>0$.
By Theorem~\ref{thm.bij}, for every $n\in \mathbf{N}$, there
exists $P^{\sigma_n}\in\mathcal{G}_0$ such that
$u_n=\psi(P^{\sigma_n})$. 
These $\sigma_n$'s are in $\Sigma_0$, and it was seen in
\cite[Lemma~1.4 and Corollary after Lemma~2.1]{marchenko}
that $\mbox{\rm supp}\,\sigma_n
       \subset[-\sqrt{\lambda},\sqrt{\lambda}]$ and
$\sigma_n(\mathbf{R})\le\lambda$ for any 
$n\in \mathbf{N}$.
\footnote{Our $\sigma_n$ and Marchenko's are different.
Namely, $\sigma_n(A)=\sigma_n^{\prime}(-A)$,
$\sigma_n^\prime$ being Marchenko's.}
Then, choosing a subsequence if necessary, we may assume
that $\sigma_n$ converges to some $\sigma\in\Sigma$ vaguely.
By Theorem~\ref{thm.conv.gen} (ii), we see that
$\psi(P^{\sigma_n})$ converges to $\psi(P^\sigma)$ 
uniformly on any bounded interval in $[0,\infty)$.
Thus $u=\psi(P^\sigma)$ on $[0,\infty)$.
\end{proof}


\end{document}